\documentclass[12pt, letterpaper, final]{amsart}

\usepackage[english]{babel}
\usepackage[latin1]{inputenc}
\usepackage{indentfirst}
\usepackage{amssymb}
\usepackage{amsmath}
\usepackage{amsthm}

\theoremstyle{definition}

\newtheorem{definicao}{{\sc Definition}}[section]

\theoremstyle{plain}

\newtheorem{proposicao}{{\sc Proposition}}[section]
\newtheorem{teorema}{{\sc Theorem}}[section]
\newtheorem{corolario}{{\sc Corollary}}[section]

\theoremstyle{remark}

\newcommand{\R}[1]{\mathbb{R}^{#1}}

\newcommand{\id}{\mathrm{id}}

\newcommand{\PI}[1]{\langle #1 \rangle}

\newcommand{\D}{\mathrm{d}}
\newcommand{\x}{\mathbf{x}}
\newcommand{\flecha}{\longrightarrow}
\newcommand{\Isom}{\mathrm{Isom}}
\newcommand{\esf}{\mathbb{S}^3}

\numberwithin{equation}{section}

\newif\ifpdf
\ifx\pdfoutput\undefined
  \pdffalse
\else {
  \pdfoutput=1
  \pdftrue
  \pdfinfo {
    /Title (Homogeneous Surfaces in S^3)
    /Creator (LaTeX)
    /Producer (LaTeX)
    /Author (Armando J. Maccori and Jose A. Verderesi)
    /Subject (Differential Geometry)
  }
}
\fi 

\ifpdf
  \usepackage[pdftex]{graphicx}
\else 
  \usepackage{graphicx}
\fi

\begin{document}

\title{Homogeneous Surfaces in $\esf$}
\author{Armando J. Maccori}
\address{
  Departamento de Matem\'atica, \newline
  \indent Universidade de S\~ao Paulo, Brazil
}
\email{maccori@ime.usp.br}
\author{Jose A. Verderesi}
\address{
  Departamento de Matem\'atica, \newline
  \indent Universidade de S\~ao Paulo, Brazil
}
\email{javerd@ime.usp.br}

\date{\today}



\begin{abstract}
The goal of this paper is to establish the classification of all
homogeneous surfaces of $3$-sphere by using the moving frame
method. We will show that such surfaces are $2$-spheres and flat
torus. 
\end{abstract}
\maketitle

\tableofcontents

\section*{Introduction}

In this work we established the classification of all homogeneous
surfaces of $\esf$ by using the method of moving frames. We will
denote the $3$-sphere by $\esf$, and it is the following subset of
$\R{4}$: 
$$ 
\esf = \big\{x \in \R{4} \mid \PI{x,x} = 1\big\}. 
$$
We say that a Riemannian surface $S$ is \emph{homogeneous} if the
group $\Isom(S)$ of all isometries of $S$ acts transitively over $S$,
i.e., if $x, y \in S$ are two distinct points of $S$, then there
exists an element $g \in \Isom(S)$ such that $y = g \cdot x$. On the
other hand, a surface $S \subset \esf$ is said to be \emph{extrinsic
homogeneous} if the group 
$$
G = \big\{g \mid g \in \Isom(\esf) \;\text{and}\; g(S) \subset S\big\} 
$$
acts transitively over $S$. Note that, if we compare these two
definitions of homogeneity, it's clear that extrinsic homogeneity
implies homogeneity. 

We obtained a Classification Theorem for immersed homogeneous surfaces
in $\esf$ (see Section \ref{cap3}). There are only two families of
homogeneous surfaces in $\esf$: the first is composed by
\emph{$2$-spheres}, given by 
$$ 
S = \big\{(x^1,x^2,x^3,x^4) \in \R{4} \mid x^4 = k, \;
(x^1)^2+(x^2)^2+(x^3)^2=1-k^2\big\}, 
$$
where $0 \le k <1$. The second family is composed by \emph{flat
  torus}. Such surfaces are given by 
$$ 
S = \big\{(x^1,x^2,x^3,x^4) \in \R{4} \mid (x^1)^2+(x^2)^2=a^2,\;
(x^3)^2+(x^4)^2=b^2\big\}, 
$$ 
where $a^2 + b^2 = 1$, and $a, b \in \R{}$. 

\section{Structure Equations of $\esf$}

Let $(e_1,e_2,e_3,e_4)$ be a moving frame of $\R{4}$ adapted to the
sphere $\esf$, i.e., $(e_1,e_2,e_3)$ belongs to $T\esf$ and $e_4(x)
= -x$. It is easy to see that $(\D e_4)_x=-\id$, since we have 
\begin{equation} \label{ses3:de4}
\D e_4 = -(\theta^1 e_1 + \theta^2 e_2 + \theta^3 e_3). 
\end{equation}

The set $(e_1,e_2,e_3,e_4)$ is an orthonormal frame, so it follows
that $\D e_i = \omega_i^k e_k$, where $\omega_i^k$ are the connection
forms of $\R{4}$. Let $(\theta^1, \theta^2, \theta^3, \theta^4)$ be
the coframe associated to $(e_1,e_2,e_3,e_4)$, i.e., $\theta^i(e_j) =
\delta_j^i$, for $i,j = 1, \ldots, 4$. 

The \emph{first structural equations} of $\R{4}$ are $\D \theta^i +
\omega_k^i \land \theta^k=0$, moreover $\D\theta^4 = 0$ over $\esf$,
and hence our set of equations reduces to 
\begin{equation} \label{1ee}
\begin{aligned}
  \D \theta^1+\omega_2^1 \land \theta^2+\omega_3^1 \land \theta^3 &=0, \\
  \D \theta^2+\omega_1^2 \land \theta^1+\omega_3^2 \land \theta^3 &=0, \\
  \D \theta^3+\omega_1^3 \land \theta^1+\omega_2^3 \land \theta^2 &=0, 
\end{aligned} 
\end{equation} 
and these equations are called the \textbf{\textit{first structural
equations}} of $\esf$. Note that $\D \theta^4 = 0$ implies the
important additional condition: 
$$ 
\omega_1^4 \land \theta^1 + \omega_2^4 \land \theta^2 + \omega_3^4
\land \theta^3 = 0. 
$$

The \emph{second structural equations} of $\R{4}$ are given by $\D
\omega_j^i + \omega_k^i \land \omega_j^k=0$. By the condition
(\ref{ses3:de4}), we obtain $\omega_4^1 = -\theta^1$, $\omega_4^2 =
-\theta^2$, $\omega_4^3 = -\theta^3$, and hence 
\begin{equation} \label{2ee}
\begin{aligned} 
  \D \omega_2^1+\omega_3^1 \land \omega_2^3 &= \theta^1 \land
  \theta^2, \\
  \D \omega_3^2+\omega_1^2 \land \omega_3^1 &= \theta^2 \land
  \theta^3, \\
  \D \omega_1^3+\omega_2^3 \land \omega_1^2 &= \theta^3 \land
  \theta^1, 
\end{aligned}
\end{equation}
and these equations are called the \textbf{\textit{second structural
equations}} of $\esf$. Moreover, the differential $2$-forms
$$ 
\begin{aligned}
\Omega_2^1 &= \theta^1 \land \theta^2, \\
\Omega_3^2 &= \theta^2 \land \theta^3, \\
\Omega_1^3 &= \theta^3 \land \theta^1,
\end{aligned}
$$
are called \emph{curvature forms} of $\esf$. 

\section{Surfaces in $\esf$} \label{cap2}

Let $S \subset \esf$ be a regular, connected, and oriented
surface. Let $(e_1,e_2,e_3)$ be an adapted orthonormal frame to $S$,
i.e., $\PI{e_i,e_j} = \delta_{ij}$, where $(e_1,e_2)$ belongs to $TS$
and $e_3 \perp TS$. We have also $\D e_i = \omega_i^k e_k$, where
$\omega_i^j$ are the connection $1$-forms of $\esf$. 

Let $(\theta^1,\theta^2,\theta^3)$ be the coframe associated to
$(e_1,e_2,e_3)$. We know that $\theta^3=0$ on $S$, because $e_3
\perp TS$, thus $\D \theta^3=0$ on $S$, and hence the equations in
\eqref{1ee} reduce to 
$$ 
\begin{aligned}
  &\D \theta^1+\omega_2^1 \land \theta^2 = 0, \\
  &\D \theta^2+\omega_1^2 \land \theta^1 = 0, \\
  &\omega_1^3 \land \theta^1+\omega_2^3 \land \theta^2 = 0,
\end{aligned} 
$$ 
and these equations are known as the \textbf{\textit{first structural
equations}} of $S$. Finally the equations in \eqref{2ee} reduce to 
$$ 
\begin{aligned} 
  \D \omega_2^1+\omega_3^1 \land \omega_2^3 &= \theta^1 \land
  \theta^2, \\ 
  \D \omega_3^2+\omega_1^2 \land \omega_3^1 &= 0, \\
  \D \omega_1^3+\omega_2^3 \land \omega_1^2 &= 0, 
\end{aligned} 
$$ 
and these are called the \textbf{\textit{second structural equations}}
of $S$. 

Now write 
$$
\begin{aligned}
\omega_1^3 &= h_{11}\theta^1 + h_{12}\theta^2, \\
\omega_2^3 &= h_{21}\theta^1 + h_{22}\theta^2,
\end{aligned}
$$
and keeping in mind the fact that 
$$ 
\omega_1^3 \land \theta^1 + \omega_2^3 \land \theta^2 = 0,
$$
it follows, by Cartan's Lemma, (see do Carmo \cite{Carmo}, page 80)
that $(h_{ij})$ is a symmetric matrix. The second fundamental form of
surface $S$ is 
$$ 
\Pi = \omega_1^3\cdot\theta^1 + \omega_2^3 \cdot \theta^2 = h_{11}
(\theta^1)^2 + 2h_{12} \theta^1\theta^2 + h_{22} (\theta^2)^2, 
$$ 
and it is, of course, a diagonalizable operator. In diagonal form,
$\omega_1^3$ and $\omega_2^3$ are written 
\begin{equation} \label{aaa}
  \begin{aligned}
    \omega_1^3 &= \lambda_1\theta^1, \\
    \omega_2^3 &= \lambda_2\theta^2,
  \end{aligned}
\end{equation} 
where $\lambda_1, \lambda_2$ are called the \emph{principal
curvatures} of $S$. 

From the second structural equations of $S$ we have 
$$
\D \omega_2^1+\omega_3^1 \land \omega_2^3 = \theta^1 \land \theta^2, 
$$
which implies 
\begin{equation} \label{gauss} 
  \D \omega_2^1 = (1+\lambda_1\lambda_2)\theta^1 \land \theta^2, 
\end{equation} 
and this is called the \emph{Gauss equation} of $S$. The function 
$$ 
K = 1 + \lambda_1 \lambda_2, 
$$ 
is the \emph{Gaussian curvature} of $S$. From another pair of
equations we have 
\begin{equation} \label{codazzi}
  \begin{aligned} 
    \D \omega_3^2-\lambda_1\omega_1^2 \land \theta^1 &=0, \\
    \D \omega_1^3-\lambda_2\omega_1^2 \land \theta^2 &=0,
  \end{aligned}
\end{equation}
called the \emph{Mainardi-Codazzi equations} of $S$. 

Differentiating the equations (\ref{aaa}), we obtain 
$$ 
\begin{aligned} 
  \D \omega_1^3=\D \lambda_1 \land \theta^1+\lambda_1 \D \theta^1, \\
  \D \omega_2^3=\D \lambda_2 \land \theta^2+\lambda_2 \D \theta^2, 
\end{aligned} 
$$ 
Therefore, by the Mainardi-Codazzi equations (\ref{codazzi}), we
conclude that 
\begin{equation} \label{codazzi2}
  \begin{aligned} 
    \lambda_2\omega_1^2 \land \theta^2 &= \D \lambda_1 \land
    \theta^1+\lambda_1 \D \theta^1, \\  
    -\lambda_1\omega_1^2 \land \theta^1 &= \D \lambda_2 \land
    \theta^2+\lambda_2 \D \theta^2,   
  \end{aligned} 
\end{equation} 
On the other hand, the first structural equations said 
\begin{equation} \label{pee}
  \begin{aligned}
    \D \theta^1  &=  \omega_1^2 \land \theta^2, \\
    \D \theta^2  &=  -\omega_1^2 \land \theta^1, \\
  \end{aligned} 
\end{equation} 
so, from (\ref{codazzi2}) and (\ref{pee}), results that the
Mainardi-Codazzi equations will be written in the form 
\begin{equation} \label{codazzi3}
  \begin{aligned}
\D \lambda_1 \land \theta^1 + (\lambda_1 - \lambda_2) \omega_1^2 \land
\theta^2 &= 0, \\ 
\D \lambda_2 \land \theta^2 + (\lambda_1 - \lambda_2) \omega_1^2 \land
\theta^1 &= 0.  
  \end{aligned}
\end{equation} 

\begin{definicao}
Let $S$ and $\tilde S$ be two surfaces of $\esf$. An \emph{isometry}
is a diffeomorphism $f:S \flecha \tilde S$ which satisfy $\PI{f_*(X),
f_*(Y)} = \PI{X, Y}$, for all pairs $X, Y \in TS$. 
\end{definicao}

\begin{proposicao} \label{ss3_prop}
Let $(e_1,\ldots,e_n)$ be a moving frame of a differentiable manifold
$M$ and let $(\theta^1,\ldots,\theta^n)$ be the coframe
associated. Then there exists a unique $1$-forms $\omega_j^i$ such
that 
\begin{equation} \label{ss3_eq1_prop}
\D \theta^i = \sum_{j=1}^n \theta^j \land \omega_j^i, 
\end{equation}
with the property $\omega_i^j = -\omega_j^i$. 
\end{proposicao}

\begin{proof}
Let $\tilde\omega_j^i$ be $1$-forms satisfying equation
\eqref{ss3_eq1_prop}. If $\omega_j^i$ also satisfies 
\eqref{ss3_eq1_prop}, then 
\begin{equation} \label{ss3_eq2_prop}
\sum_{j=1}^n \theta^j \land (\omega_j^i - \tilde\omega_j^i) = 0. 
\end{equation}
By Cartan's Lemma, from equation \eqref{ss3_eq2_prop}, follows that 
\begin{equation} \label{ss3_eq3_prop}
\omega_j^i - \tilde\omega_j^i = \sum_{k=1}^n a_{jk}^i \theta^k, 
\end{equation}
where $a_{jk}^i$ are symmetric ($a_{jk}^i = a_{kj}^i$). Since
$(\theta^1, \ldots, \theta^n)$ is a base for the set of $1$-forms in
$M$, then there exist $\Gamma_{jk}^i$ such that 
\begin{equation} \label{ss3_eq4_prop}
\tilde\omega_j^i = \sum_{k=1}^n \Gamma_{jk}^i \theta^k. 
\end{equation}
Thus, from equations \eqref{ss3_eq3_prop} and \eqref{ss3_eq4_prop},
follows that 
\begin{equation} \label{ss3_eq5_prop}
\omega_j^i = \sum_{k=1}^n (\Gamma_{jk}^i + a_{jk}^i)\theta^k. 
\end{equation}
If $\omega_i^j = -\omega_j^i$ is satisfied, then the equation 
\eqref{ss3_eq5_prop} assumes the form 
$$ 
(\Gamma_{jk}^i + a_{jk}^i) + (\Gamma_{ik}^j + a_{ik}^j) = 0,
$$ 
which is equivalent to 
\begin{equation} \label{ss3_eq6_prop}
a_{jk}^i + a_{ik}^j = -(\Gamma_{jk}^i + \Gamma_{ik}^j). 
\end{equation}
Cyclic permuting the indices $i,j,k$ in equation \eqref{ss3_eq6_prop},
we write 
\begin{equation} \label{ss3_eq7_prop}
\begin{aligned}
  a_{jk}^i + a_{ik}^j &= -(\Gamma_{jk}^i + \Gamma_{ik}^j), \\
  a_{ij}^k + a_{kj}^i &= -(\Gamma_{ij}^k + \Gamma_{kj}^i), \\
  a_{ki}^j + a_{ji}^k &= -(\Gamma_{ki}^j + \Gamma_{ji}^k). 
\end{aligned}
\end{equation}
In \eqref{ss3_eq7_prop}, if we add the first equation with the second
and subtract from the third equation, in both members, we will obtain
(considering the fact that $a_{jk}^i$ are symmetric) 
$$ 
a_{jk}^i = \frac{1}{2}(\Gamma_{ki}^j + \Gamma_{ji}^k - \Gamma_{ij}^k -
\Gamma_{kj}^i - \Gamma_{jk}^i -\Gamma_{ik}^j). 
$$ 
It follows that 
$$ 
\omega_j^i = \frac{1}{2} \sum_{k=1}^n (\Gamma_{jk}^i + \Gamma_{ki}^j +
\Gamma_{ji}^k - \Gamma_{kj}^i - \Gamma_{ik}^j - \Gamma_{ij}^k)
\theta^k, 
$$ 
and note that $\omega_i^j = -\omega_j^i$. This demonstrates the
existence and uniqueness of connection $1$-forms $\omega_j^i$. 
\end{proof}

\begin{corolario} \label{ss3_cor}
Let $M$ and $\tilde M$ be two Riemannian manifolds of dimension $n$
and let $f:M \flecha \tilde M$ be an isometry. Let $(\tilde\theta^1,
\ldots, \tilde\theta^n)$ be an adapted coframe in $\tilde M$ whose
connection forms are $\tilde\omega_j^i$. If $\omega_j^i$ are the
connection forms of $M$ in the adapted coframe $(\theta^1, \ldots,
\theta^n)$ where $\theta^i = f^*\tilde\theta^i$, for $i = 1, \ldots,
n$, then $\omega_j^i = f^*\tilde\omega_j^i$. 
\end{corolario}

\begin{proof}
According to Proposition \ref{ss3_prop}, in $\tilde M$, the connection
$1$-forms $\tilde\omega_j^i$ are the only one satisfying the
structural equations 
\begin{equation} \label{ss3_eq1_cor}
\D \tilde\theta^i + \sum_{j=1}^n \tilde\omega_j^i \land
\tilde\theta^i = 0. 
\end{equation}
Again, by Proposition \ref{ss3_prop} applied to $M$, the connection
$1$-forms $\omega_j^i$ are the only one satisfying the structural
equations 
\begin{equation} \label{ss3_eq2_cor}
\D \theta^i + \sum_{j=1}^n \omega_j^i \land \theta^i = 0, 
\end{equation}
where $\theta^i = f^*\tilde\theta^i$, for $i=1, \ldots, n$. 

Calculating the pullback $f^*$ in equation \eqref{ss3_eq1_cor} and
comparing with \eqref{ss3_eq2_cor}, by uniqueness of connection forms,
we conclude that 
$$ 
\omega_j^i = f^*\tilde\omega_j^i, 
$$ 
as we wished. 
\end{proof}

\begin{teorema}
Let $S$ and $\tilde S$ be two surfaces in $\esf$ and let $f:\esf
\flecha \esf$ be an isometry such that $f(S) \subset \tilde S$. In these
conditions, the following assertions are true: 
\begin{enumerate}
\item[(i)] If $K$ and $\tilde K$ are the Gaussian curvatures of $S$
  and $\tilde S$, respectively, then $K(p) = \tilde K\big(f(p)\big)$,
  for all $p \in S$. 
\item[(ii)] If $\lambda_1, \lambda_2$ and $\tilde\lambda_1,
  \tilde\lambda_2$ are the principal curvatures of $S$ and $\tilde S$,
  respectively, then $\lambda_1(p) =  \tilde\lambda_1\big(f(p)\big)$
  and $\lambda_2(p) = \tilde\lambda_2\big(f(p)\big)$, for all $p \in
  S$. 
\end{enumerate}
\end{teorema}

\begin{proof}
Let us prove (i). Let $(\tilde\theta^1, \tilde\theta^2,
\tilde\theta^3)$ be an adapted orthonormal coframe to the surface
$\tilde S$. From the first structural equations, we detach the
following 
\begin{equation} \label{sh1a}
\begin{aligned}
\D \tilde\theta^1 + \tilde\omega_2^1 \land \tilde\theta^2 &= 0, \\
\D \tilde\theta^2 + \tilde\omega_1^2 \land \tilde\theta^1 &= 0,
\end{aligned}
\end{equation} 
and
\begin{equation} \label{sh1b}
\D \tilde\omega_2^1 = \tilde K \,\tilde\theta^1 \land \tilde\theta^2. 
\end{equation} 
Since $f$ is an isometry, the ternary $(\theta^1,\theta^2,\theta^3)$,
where $\theta^i = f^*(\tilde\theta^i)$, for $i = 1,2,3$, defines an
adapted orthonormal coframe to the surface $S$. Again, for this
coframe, we detach the following structural equations 
\begin{equation} \label{sh2a}
\begin{aligned}
\D \theta^1 + \omega_2^1 \land \theta^2 &= 0, \\
\D \theta^2 + \omega_1^2 \land \theta^1 &= 0, 
\end{aligned}
\end{equation}
and
\begin{equation} \label{sh2b}
\D \omega_2^1 =  K \,\theta^1 \land \theta^2.
\end{equation}
Applying the pullback $f^*$ in equations \eqref{sh1a} and
\eqref{sh1b}, we obtain 
\begin{equation} \label{sh3a}
\begin{aligned}
\D \theta^1 + f^*(\tilde\omega_2^1) \land \theta^2 &= 0, \\
\D \theta^2 + f^*(\tilde\omega_1^2) \land \theta^1 &= 0, \\
\end{aligned}
\end{equation}
and
\begin{equation} \label{sh3b}
\D \big(f^*(\tilde\omega_2^1)\big) = \tilde K(f) \,\theta^1 \land
\theta^2. 
\end{equation}
By Proposition \ref{ss3_prop} and according to the first structural
equations in \eqref{sh2a} and \eqref{sh3a}, we conclude that
$\omega_2^1 = f^*(\tilde\omega_2^1)$, which together with the Gauss
equations \eqref{sh2b} and \eqref{sh3b}, generates $K = \tilde K(f)$,
i.e., $K(p) = \tilde K\big( f(p) \big)$, for all $p \in S$. 

Finally, to prove item (ii), let $(\tilde\theta^1, \tilde\theta^2,
\tilde\theta^3)$ be an adapted orthonormal coframe to the surface
$\tilde S$. In an analogous way, the set $\theta^i = f^*(\tilde
\theta^i)$ for $i = 1,2,3$ is an adapted orthonormal coframe to the
surface $S$. We have the following expressions 
\begin{equation} \label{sh4}
\D \tilde\theta^3 + \tilde\omega_1^3 \land \tilde\theta^1 +
\tilde\omega_2^3 \land \tilde\theta^2 = 0, 
\end{equation}
and  
\begin{equation} \label{sh5}
\D \theta^3 + \omega_1^3 \land \theta^1 + \omega_2^3 \land \theta^2 =
0. 
\end{equation}
By Corollary \ref{ss3_cor}, and applying the pullback $f^*$ in
equation \eqref{sh4}, we obtain, by direct comparison with equation
\eqref{sh5}, 
\begin{equation} \label{sh6}
\omega_1^3 = f^* (\tilde\omega_1^3) 
\;\;\;\text{and}\;\;\; 
\omega_2^3 = f^* (\tilde\omega_2^3). 
\end{equation} 
Now, writing 
\begin{equation} \label{sh7}
\omega_1^3 = \lambda_1 \theta^1 
\;\;\;\text{and}\;\;\;
\omega_2^3 = \lambda_2 \theta^2,
\end{equation} 
while we also have 
\begin{equation} \label{sh8}
\tilde\omega_1^3 = \tilde\lambda_1 \tilde\theta^1
\;\;\;\text{and}\;\;\;
\tilde\omega_2^3 = \tilde\lambda_2 \tilde\theta^2.
\end{equation}
Applying $f^*$ in \eqref{sh8} and using \eqref{sh6}, we conclude 
\begin{equation} \label{sh9}
\omega_1^3 = \tilde\lambda_1(f)\, \theta^1 
\;\;\;\text{and}\;\;\;
\omega_2^3 = \tilde\lambda_2(f)\, \theta^2,
\end{equation}
and therefore, comparing the last equations with expressions in
\eqref{sh7}, we determine that 
$$ 
\lambda_1 = \tilde\lambda_1(f) 
\;\;\;\text{and}\;\;\; 
\lambda_2 = \tilde\lambda_2(f),
$$ 
this immediately implies that $\lambda_1(p) =
\tilde\lambda_1\big(f(p)\big) $ and $\lambda_2(p) =
\tilde\lambda_2\big(f(p)\big)$, for all $p \in S$. 
\end{proof}

\begin{corolario} \label{hcpc}
If $S$ is an extrinsic homogeneous surface of $\esf$, then we have 
\begin{enumerate}
\item[(i)] its Gaussian curvature $K$ is a constant. 
\item[(ii)] its principal curvatures $\lambda_1$ and $\lambda_2$ are
  constant functions. 
\end{enumerate}
\end{corolario}

\begin{proof}
Item (i). Since $S$ is a homogeneous surface, then for any $p, q \in
S$ there exists an isometry $f:\esf \flecha \esf$ such that $f(p) =
q$. Now, by the last Theorem, putting $\tilde S = S$ (and hence
$\tilde K = K$), we have $K(p) = K\big(f(p)\big)$, for all $p \in
S$. Therefore, it results that $K(p) = K(q)$, for all $p, q \in S$ and
thus $K(p)$ is constant on $S$. 

Item (ii). In a similar way of the item (i), letting $\tilde S = S$
and taking $p, q \in S$, there exists $f \in \Isom(\esf)$ such that $q
= f(p)$. Thus we will have $\lambda_1(p) = \lambda_1\big(f(p)\big) =
\lambda_1(q)$ and $\lambda_2(p) = \lambda_2\big(f(p)\big) =
\lambda_2(q)$, therefore $\lambda_1$ and $\lambda_2$ are constants on
$S$. 
\end{proof}

\section{Classification of Homogeneous Surfaces of $\esf$}
\label{cap3} 

\begin{proposicao}
If $S$ is an umbilic surface, then its principal curvatures are
constants. In particular, $K$ will be constant. 
\end{proposicao} 
\begin{proof}
Since $S$ is umbilic, i.e., $\lambda_1 = \lambda_2 = \lambda$, the
equations in (\ref{codazzi3}) reduce to $\D \lambda \land \theta^j =
0$, for $j = 1, 2$. Thus $\D \lambda = 0$, and hence $\lambda$ is
constant. Moreover, since $K = 1+\lambda^2$, it results that $K$ is
also a constant. 
\end{proof}

\begin{proposicao} 
If the principal curvatures of $S$ are constants, then or $S$ is
umbilic or $S$ has null Gaussian curvature $K$. 
\end{proposicao} 
\begin{proof}
Suppose that $\lambda_1$ and $\lambda_2$ are constants. From the
equations in (\ref{codazzi3}) it follows that $(\lambda_1 - \lambda_2)
\omega_1^2 \land \theta^j = 0$, for $j = 1, 2$. Thus, there are two
possibilities: or $\lambda_1=\lambda_2$, and hence $S$ is umbilic, or
$\omega_1^2 \land \theta^j = 0$, for $j = 1, 2$, which implies
$\omega_1^2=0$, and hence $\D \omega_1^2 = 0$. However, looking to the
equation (\ref{gauss}) we conclude that $K\theta^1 \land \theta^2 =
0$. Therefore $K = 0$, and hence $S$ is a surface with null Gaussian
curvature. 
\end{proof} 

\subsection{The case $\lambda_1=\lambda_2$}

Let $S \subset \esf$ be a homogeneous surface, and hence, by Corollary
\ref{hcpc}, its principal curvatures $\lambda_1 = \lambda_2 = \lambda
\in \R{}$. Let $(e_1, e_2)$ be the principal direction of $S$ and
$e_3$ the normal field of $S$. Remember that in these conditions,
$\theta^3 \equiv 0$ over $S$. Since the frame $(e_1,e_2,e_3)$ is
adapted to $S$, we have the following set of equations 
$$ 
\begin{aligned} 
  \D e_1 &= \omega_1^2 e_2-\lambda\theta^1e_3+\theta^1e_4, \\ 
  \D e_2 &= \omega_2^1e_1-\lambda\theta^2e_3+\theta^2e_4, \\ 
  \D e_3 &= \lambda\theta^1e_1+\lambda\theta^2e_2=\lambda\,\id, \\ 
  \D e_4 &= -\theta^1e_1-\theta^2e_2 = -\id,
\end{aligned} 
$$ 
moreover, note that $\omega_3^1 = \lambda\theta^1$ and $\omega_3^2 =
\lambda\theta^2 $. 

Then, consider the following vector field defined on $S$ 
$$ 
X = \x -\frac{1}{\lambda}e_3, 
$$ 
where $\x$ is a parametrization of $S$. We will show that $X = \x_0 =
\mathrm{constant}$ on $S$. In fact, differentiating $X$ we obtain 
\begin{align*} 
  \D X & = \D \Big{(}\x-\frac{1}{\lambda}e_3\Big{)} \\
  & = \theta^1e_1 + \theta^2e_2 - \frac{1}{\lambda}(\lambda\theta^1e_1
  + \lambda\theta^2e_2) \\
  & = \theta^1e_1+\theta^2e_2-\theta^1e_1-\theta^2e_2 = 0. 
\end{align*} 
Thus, $\D X = 0$ on $S$, and hence $X = \x_0 = \mathrm{constant}$. If
we write 
$$ 
\x_0 = \x-\frac{1}{\lambda}e_3 \quad \textnormal{which implies} \quad
\x-\x_0 = \frac{1}{\lambda}e_3, 
$$  
Taking the norm, in both members, on the last equation, we conclude 
\begin{equation} \label{e1} 
  ||\x-\x_0||=\frac{1}{\lambda}, 
\end{equation}  
which immediately implies that 
\begin{equation} \label{e2} 
  \PI{\x-\x_0,\x-\x_0}=\frac{1}{\lambda^2}, 
\end{equation} 
and this is an equation of a sphere with center $\x_0$ and radius
$\frac{1}{|\lambda|}$. 

Remember that $S$ is a connected homogeneous surface, therefore it is
a \emph{complete} surface. Hence $S$ is a whole $2$-sphere. 

\subsection{The case $\lambda_1 \ne \lambda_2$}

Let $S \subset \esf$ be a surface, and suppose its principal
curvatures $\lambda_1, \lambda_2$ are constant and distinct. From
equations of Mainardi-Codazzi (\ref{codazzi3}), we have 
$$ 
\begin{aligned} 
  (\lambda_1-\lambda_2)\omega_1^2 \land \theta^2 &=0, \\ 
  (\lambda_1-\lambda_2)\omega_1^2 \land \theta^1 &=0.
\end{aligned}
$$ 
Therefore, it follows that $\omega_1^2=0$. From equation of Gauss
(\ref{gauss}), follows that $K\theta^1 \land \theta^2 = 0$, which
implies that $K = 0$, but since $K = 1 + \lambda_1\lambda_2$, it
results the condition: 
$$
\lambda_1\lambda_2 = -1. 
$$ 
Since $K = 0$, there exists a parametrization $\x: U \ni (u^1, u^2)
\longmapsto \R{4}$, where $U$ is an open set of $\R{2}$, $\x(U)
\subset S$, and satisfying the property 
$$ 
\frac{\partial \x }{\partial u^1} = e_1 
\quad \textnormal{and} \quad
\frac{\partial \x }{\partial u^2}=e_2. 
$$
In these conditions, we know that 
\begin{equation}
\D u^1 = \theta^1 
\quad \textnormal{and} \quad 
\D u^2 = \theta^2, 
\end{equation}
since $(\theta^1,\theta^2,\theta^3,\theta^4)$ is the dual base of
$(e_1,e_2,e_3,e_4)$. 

On the other hand, we have the following system of equations 
$$ 
\begin{aligned} 
  \D e_1 &= \lambda_1\theta^1e_3+\theta^1e_4 = (\lambda_1e_3 + e_4)\D
  u^1, \\ 
  \D e_2 &= \lambda_2\theta^2e_3+\theta^2e_4 = (\lambda_2e_3 + e_4)\D
  u^2, \\ 
  \D e_3 &= -\lambda_1\theta^1e_1-\lambda_2\theta^2e_2 =
  -\lambda_1 e_1\,\D u^1 - \lambda_2 e_2\,\D u^2, \\
  \D e_4 &= -\theta^1e_1-\theta^2e_2 = -e_1\,\D u^1 - e_2\,\D u^2. 
\end{aligned} 
$$ 

Then, let us consider the following vector fields 
$$ 
\begin{aligned}
  f_1 &= e_1, \\
  f_2 &= e_2, \\
  f_3 &= \frac{\lambda_1e_3+e_4}{\sqrt{\lambda_1^2+1}}, \\
  f_4 &= \frac{\lambda_2e_3+e_4}{\sqrt{\lambda_2^2+1}}, 
\end{aligned}
$$ 
it is easy to check that $\PI{f_i,f_j}=\delta_{ij}$, i.e.,
the set $(f_1,f_2,f_3,f_4)$ forms a base of $\R{4}$. 

Differentiating the vector fields $f_i$, we obtain 
$$ 
\begin{aligned} 
  \D f_1 &= \D e_1 = \sqrt{\lambda_1^2+1}\,\theta^1 f_3 =
  \Big(\sqrt{\lambda_1^2+1}\,\D u^1\Big) f_3, \\
  \D f_2 &= \D e_2 = \sqrt{\lambda_2^2+1}\,\theta^2 f_4 =
  \Big(\sqrt{\lambda_2^2+1}\,\D u^2 \Big)f_4, \\
  \D f_3 &= \D
  \bigg{(}\dfrac{\lambda_1e_3+e_4}{\sqrt{\lambda_1^2+1}}\bigg{)} =
  -\sqrt{\lambda_1^2+1}\,\theta^1f_1 = -\Big(\sqrt{\lambda_1^2+1}\,\D
  u^1\Big)f_1, \\ 
  \D f_4 &= \D
  \bigg{(}\dfrac{\lambda_2e_3+e_4}{\sqrt{\lambda_2^2+1}}\bigg{)} =
  -\sqrt{\lambda_2^2+1}\,\theta^2f_2 = -\Big(\sqrt{\lambda_2^2+1}\,\D
  u^2\Big)f_2, 
\end{aligned}
$$ 
and we will denote by $k_i = \sqrt{\lambda_i^2+1}$, for $i = 1,2$,
just for simplify the expressions above. So, we rewrite 
\begin{equation} \label{ccc} 
  \begin{aligned} 
    \D f_1 &= k_1\,\D u^1\,f_3, \\ 
    \D f_2 &= k_2\,\D u^2\,f_4, \\ 
    \D f_3 &= -k_1\,\D u^1\,f_1, \\ 
    \D f_4 &= -k_2\,\D u^2\,f_2. 
  \end{aligned} 
\end{equation} 

Thus, if we observe the equations in (\ref{ccc}) and the fact that
$e_i = f_i$ for $i=1,2$, we conclude that 
$$
\begin{array}{ll}
  \begin{cases}
    \frac{\partial f_1}{\partial u^1} &= k_1f_3, \\ 
    \frac{\partial f_1}{\partial u^2} &= 0,
  \end{cases} &
  \begin{cases}
    \frac{\partial f_2}{\partial u^1} &= 0, \\
    \frac{\partial f_2}{\partial u^2} &= k_2f_4,
  \end{cases} 
  \\ \\
  \begin{cases}
    \frac{\partial f_3}{\partial u^1} &= -k_1f_1, \\
    \frac{\partial f_3}{\partial u^2} &= 0, 
  \end{cases} & 
  \begin{cases}
    \frac{\partial f_4}{\partial u^1} &= 0. \\
    \frac{\partial f_4}{\partial u^2} &= -k_2f_2,
  \end{cases}
\end{array}
$$

Consider the following curve on $S$: 
$$ 
c_1(u^1) = \x (u^1,u_0^2), 
$$
where $u_0^2$ is fixed. Then 
$$ 
X = c_1(u^1)+\frac{1}{k_1}f_3, 
$$ 
is a vector field defined on $S$. We have now 
\begin{align*} 
  \dfrac{\partial X}{\partial u^1} & = \dfrac{\partial}{\partial
    u^1}\bigg(c_1(u^1)+\dfrac{1}{k_1}f_3\bigg) \\
  & = \dfrac{\partial \x }{\partial u^1}(u^1,u_0^2) +
  \dfrac{1}{k_1}\dfrac{\partial f_3}{\partial u^1} \\
  & = f_1 + \dfrac{1}{k_1}(-k_1f_1) = 0. 
\end{align*} 
Therefore $X$ is a constant vector field, i.e., we can write 
$$ 
c_1(u^1)+\frac{1}{k_1}f_3=p_0, 
$$ 
where $p_0$ is a point of $S$. 

It follows that 
$$ 
c_1(u^1)+\frac{1}{k_1}f_3 = p_0 
\;\;\textnormal{which implies}\;\;
c_1(u^1)-p_0 = -\frac{1}{k_1}f_3; 
$$  
and taking the norm, in both members, we obtain 
$$ 
\big\Vert c_1(u^1)-p_0\big\Vert = \frac{1}{k_1}||f_3|| 
\;\;\textnormal{which implies}\;\;
\big\Vert c_1(u^1)-p_0\big\Vert = \frac{1}{k_1}, 
$$ 
and this is an equation of a circle of center $p_0$ and radius
$\frac{1}{k_1}$. 

We will show that this circle belongs to a special plane. In fact, we
have 
\begin{align*} 
  \dfrac{\partial}{\partial u^1}\PI{c_1(u^1),f_2} & =
  \Big\langle\dfrac{\partial \x }{\partial u^1},f_2\Big\rangle + 
  \Big\langle\x (u^1,u_0^2),\dfrac{\partial f_2}{\partial u^1}\Big\rangle \\
  & = \PI{f_1,f_2}+\big\langle c_1(u^1),0\big\rangle = 0, 
\end{align*} 
and hence $\big\langle c_1(u^1), f_2\big\rangle = a$, where $a \in
\R{}$. 

On the other hand, 
\begin{align*} 
  \dfrac{\partial}{\partial u^1}\PI{c_1(u^1),f_4} & =
  \Big\langle\dfrac{\partial \x}{\partial u^1},f_4\Big\rangle +
  \Big\langle\x (u^1,u_0^2),\dfrac{\partial f_4}{\partial
  u^1}\Big\rangle \\ 
  & = \PI{f_1,f_4}+\big\langle c_1(u^1),0\big\rangle = 0, 
\end{align*} 
and hence $\big\langle c_1(u^1), f_4\big\rangle = b$, where $b \in
\R{}$. 

The last pair of equations defines a plane $\pi_1$ in $\R{4}$, i.e., 
$$ 
\pi_1: \;\; \big\langle c_1(u^1),f_2\big\rangle = a, \quad \big\langle
c_1(u^1),f_4\big\rangle = b, \qquad a, b \in \R{}. 
$$  

Therefore, it results that the circle 
$$ 
\mathcal{C}_1: \;\; \big\langle c_1(u^1)-p_0,c_1(u^1)-p_0\big\rangle =
\frac{1}{k_1^2}, 
$$ 
belongs to plane $\pi_1$. 

In a similar way, consider the following curve in $S$: 
$$ 
c_2(u^2)=\x (u_0^1,u^2), 
$$
where $u_0^1$ is fixed. Let 
$$ 
Y=c_2(u^2)+\frac{1}{k_2}f_4, 
$$ 
be a vector field on $S$. We have 
\begin{align*} 
  \dfrac{\partial Y}{\partial u^2} & = \dfrac{\partial}{\partial
    u^2}\bigg(c_2(u^2)+\dfrac{1}{k_2}f_4\bigg) \\
  & = \dfrac{\partial \x }{\partial
    u^2}(u_0^1,u^2) + \dfrac{1}{k_2}\dfrac{\partial f_4}{\partial u^2}
  \\
  & = f_2+\dfrac{1}{k_2}(-k_2f_2) = 0. 
\end{align*} 
Therefore $Y$ is a constant vector field, i.e., we can write 
$$ 
c_2(u^2)+\frac{1}{k_2}f_4=q_0, 
$$ 
where $q_0$ is a point of $S$. 

It follows that, 
$$ 
c_2(u^2) + \frac{1}{k_2}f_4 = p_0 
\;\;\textnormal{which implies}\;\; 
c_2(u^2) - q_0 = -\frac{1}{k_2}f_4; 
$$ 
and again, taking the norm, in both members, we obtain 
$$ 
\big\Vert c_2(u^2)-q_0\big\Vert = \frac{1}{k_2}||f_4|| 
\;\;\textnormal{which implies}\;\;
\big\Vert c_2(u^2)-q_0\big\Vert = \frac{1}{k_2}, 
$$  
and this is an equation of a circle centered in $q_0$ and with radius
$\frac{1}{k_2}$. 

Again, we will show that this circle belongs to a special plane. In
fact, we have 
\begin{align*} 
  \dfrac{\partial}{\partial u^2}\PI{c_2(u^2),f_1} & = 
  \Big\langle\dfrac{\partial \x }{\partial u^2},f_1\Big\rangle + \Big\langle\x
  (u_0^1,u^2),\dfrac{\partial f_1}{\partial u^2}\Big\rangle  \\ 
  & = \PI{f_2,f_1}+\big\langle c_2(u^2),0\big\rangle = 0, 
\end{align*}
and hence $\big\langle c_2(u^2),f_1\big\rangle=c$, where $c \in
\R{}$. 

On the other hand, 
\begin{align*} 
  \dfrac{\partial}{\partial u^2}\PI{c_2(u^2),f_3} & =
  \Big\langle\dfrac{\partial \x}{\partial u^2},f_3\Big\rangle + 
  \Big\langle\x (u_0^1,u^2),\dfrac{\partial f_3}{\partial u^2}\Big\rangle \\ 
  & = \PI{f_2,f_1}+\big\langle c_2(u^2),0\big\rangle = 0, 
\end{align*}
and hence $\big\langle c_2(u^2),f_3\big\rangle=d$, where $d \in
\R{}$. 

The last pair of equations defines a plane $\pi_2$ in $\R{4}$, i.e., 
$$ 
\pi_2: \;\; \big\langle c_2(u^2),f_1\big\rangle=c, \quad \big\langle
c_2(u^2),f_3\big\rangle = d, \qquad c,d \in \R{}. 
$$  

Finally, it results that the circle 
$$ 
\mathcal{C}_2: \;\; \big\langle c_2(u^2)-q_0,c_2(u^2)-q_0\big\rangle =
\frac{1}{k_2^2}, 
$$ 
belongs to the plane $\pi_2$. 

Since $\{f_1,f_3\}$ and $\{f_2,f_4\}$ belong to mutually orthogonal
planes, the circles $\mathcal{C}_1$ and $\mathcal{C}_2$ are both
orthogonal, and hence they generate a torus in $\esf$. 

Another time, keeping in mind that $S$ is a connected homogeneous
surface, it follows that $S$ is a \emph{complete} surface, and hence
$S$ is a whole torus. 

\section{Conclusions}

From Corollary \ref{hcpc} and from both cases (i), (ii), we conclude
that complete immersed surfaces with constant principal curvatures are
$2$-spheres and torus. Since these surfaces are homogeneous, we have
the following classification theorem. 

\begin{teorema}
If $S$ is a regular connected immersed homogeneous surface of $\esf$,
then $S$ is one and only one surface between the following types: 
\begin{itemize}
\item[(i)] $S$ is an immersed $2$-sphere in $\esf$.
\item[(ii)] $S$ is an immersed flat torus $(\cong$ $S^1 \times S^1)$
  in $\esf$. 
\end{itemize}
\end{teorema}



\end{document}